\documentclass[11pt, leqno]{article}
\usepackage{pdfsync}
\usepackage{array} 
\usepackage{amsfonts}
\usepackage{amsmath}
\usepackage{amssymb}
\usepackage{graphicx}
\usepackage{color}
\usepackage{bm}
\usepackage[utf8]{inputenc}
\usepackage[T1]{fontenc}
\newtheorem{theorem}{Theorem}[section]
\newtheorem{corollary}[theorem]{Corollary}

\newtheorem{proposition}[theorem]{Proposition}

\newcommand{\re}{\mathbb{R}}

\newcommand{\ren}{\mathbb{R}^n}

\newcommand{\Hn}{\mathbb{H}^n}
\newcommand\C{\mathbb{C}}
\renewcommand\H{\mathbb{H}}
\renewcommand\O{\mathbb{O}}

\newcommand{\ve}{{\varepsilon}}
\newcommand{\Sh}{\mbox{sinh\,}}
\newcommand{\Ch}{\mbox{cosh\,}}
\newcommand{\Cth}{\mbox{coth\,}}
\newcommand\ssf{\hspace{.25mm}}
\newcommand\msb{\hspace{-.5mm}}

\newcommand\vsf{\hspace{.1mm}}

\def\qed{\unskip\kern 6pt \penalty 500
\raise -2pt\hbox{\vrule \vbox to8pt{\hrule width 6pt \vfill\hrule}\vrule}\par}

\numberwithin{equation}{section}

\begin{document}

\title{\textbf{Asymptotic behaviour  for the Heat \\ Equation
in Hyperbolic Space}}

\author{{\Large Juan Luis V\'azquez}\\ [4pt]
Universidad Aut\'{o}noma de Madrid\\
} 
%

\maketitle

\begin{abstract}
Following the classical result of long-time asymptotic convergence towards the Gaussian kernel that holds true for integrable solutions of the Heat Equation posed in the Euclidean Space $\mathbb{R}^n$, we examine the question of long-time behaviour of the Heat Equation in the Hyperbolic Space $\mathbb{H}^n$, $n>1$, also for integrable solutions. We show that the typical convergence proof towards the fundamental solution works in the class of radially symmetric solutions. We also prove the more precise result that says that this limit behaviour is exactly described by the 1D Euclidean kernel, but only after correction of  a remarkable outward drift with constant speed produced by the geometry. Finally, we find that such fine convergence results are false for general nonnegative solutions with integrable initial data.
\end{abstract}

\

\

\noindent \it 2010 Mathematics Subject Classification\rm: 35K05, 35K08, 35B40, 35R01.


\noindent\it Keywords and phrases\rm: Heat equation, hyperbolic space, asymptotic behaviour, equations with drift.

\newpage



\newpage



\section{Introduction}\label{sec.conv1}
We study the long-time behaviour of the solutions of the Heat Equation posed in the Hyperbolic Space $\mathbb{H}^n$, $n>1$. We  take as motivation the classical result of asymptotic convergence towards the Gaussian kernel  for the Heat Equation $\partial _t u=\Delta u$ posed in the Euclidean Space $\ren$, $n\ge 1$, which holds for all integrable initial data $u(x,0)\in L^1(\ren)$ and admits finer versions with explicit rate of convergence under extra assumptions on the data, cf. \cite{EvansPDE, VazGN}. We recall next the main versions of the Euclidean statement.

\begin{theorem} \label{thm.conv.euc1} Let  $u(x,t)$ be a solution of the heat equation  taking initial data $u_0\in L^1(\ren)$ with mass $\int_{\ren} u_0\,dx=M>0$. Then  we have
\begin{equation}
\displaystyle \|u(x,t)-MG_t(x)\|_{L^1(\ren)}\to 0 \quad \mbox{as \ $t\to\infty$}.
\end{equation}
and
\begin{equation}
\displaystyle t^{n/2}\|u(x,t)-MG_t(x)\|_{L^{\infty}(\ren)}\to 0 \quad \mbox{as \ $t\to\infty$}.
\end{equation}
Here, $G_t$ is the Gaussian kernel, $G_t(x)=(4\pi t)^{-n/2} e^{-|x|^2/4t}$. The rates are optimal in that generality. Norms are taken for fixed $t$.
\end{theorem}

This result has a stronger version under stricter conditions.

\begin{theorem} \label{thm.conv.euc2} Assume that $u_0\in L^1(\ren)$ and is compactly supported (or has finite second moment) and assume that we put the origin of coordinates at the center of mass of $u_0$. Then if  $u(x,t)$ is the solution of the heat equation we have
\begin{equation}
\displaystyle \|u(x,t)-MG_t(x)\|_{L^1(\ren)}\le C\,t^{-1} \quad \mbox{for all \ $t\ge 1$}.
\end{equation}
\end{theorem}

In this paper we examine the question of whether a similar convergence holds when the Heat Equation is posed in the Hyperbolic Space.  As a positive result, we show that  convergence towards the fundamental solution works in the class of radially symmetric solutions with finite mass.

\begin{theorem} \label{thm.conv.hypb.pos}  Let $u(x,t)$ be a solution of the heat equation in the hyperbolic space whose initial function $u_0\in L^1(\mathbb{H}^n)$ is radially symmetric in geodesic polar coordinates around a pole $o\in \Hn$ and has mass $M$. Then we have
\begin{equation}
\displaystyle \|u(x,t)-MP_t(x)\|_{L^1(\mathbb{H}^n)}\to 0 \quad \mbox{as \ $t\to\infty$}\,,
\end{equation}
where $P_t(x)$ denotes the fundamental solution (h.f.s. for short) of the heat equation in $\Hn$ centered at the pole. Moreover,
\begin{equation}
\|u(x,t)-MP(x,t)\|_{L^{\infty}(\mathbb{H}^n)}=O(t^{-3/2}e^{-\lambda_1 t})\quad
\mbox{as $t\to\infty$}.
\end{equation}
\end{theorem}

Here mass means $M=\int_{\Hn} u_0(x)\,d\mu(x)$, where $d\mu$ is the volume element in $\Hn$. Note that we are correcting the norm of the last formula by the precise size of the h.f.s., $\|P_t(\cdot)\|_\infty=O(t^{-3/2}e^{-(n-1)^2t/4})$. By interpolation, convergence estimates with the appropriate weights are obtained in all $L^p(\Hn)$ norms, for all $1<p<\infty$, see formula \eqref{rd.convp}.

For the reader's benefit, the relevant needed concepts and facts are summarized in Section \ref{sec.prelim}. In Section \ref{sec.pbhi} we prove that fundamental solutions bound any solution with compactly supported initial data from above and below on expanding sets,
and we derive a Harnack inequality of elliptic type for such solutions.

We study the peculiar behaviour of the h.f.s. for large times in Section \ref{sec.mass} based on the mass analysis, and we introduce the concepts of mass escape and mass line, which happens to be $r(t)=(n-1)t$ for large $t$. It is the realization of a drift with constant speed $c=(n-1)$ in the forward direction of the radial axis  polar coordinates. This drift had been noted by Davies \cite{Davies} and makes a big difference with the behaviour found in Euclidean space.

The  asymptotic behaviour of the h.f.s. is made precise in Section \ref{sec.1Ddrift} where we prove that the fundamental solution itself resembles for large times the fundamental solution of the heat equation in 1D (the Gaussian kernel), but for the drift. By means of a careful analysis of the equations we  prove asymptotic convergence to the 1D Euclidean profile up to this drift, see Theorems \ref{conv.fsto1d}, \ref{conv.fsnto1d}. After that work, the asymptotic convergence results for radial solutions of Theorem \ref{thm.conv.hypb.pos} are proved in Section~\ref{sec.rad}.

Section \ref{sec.horo} studies the class of horospheric solutions. They are not integrable in $\Hn$ but their analysis is simpler. They allow to derive in a very direct way the corresponding drift effect that causes the mass to escape into the outer field with constant speed.

In the last part of the paper we prove a negative result. We find that the previous kind of stabilization towards the fundamental solution does not hold for general nonnegative solutions of the heat flow with integrable initial data if radial symmetry is lacking. Here is the  result.

\begin{theorem} \label{thm.conv.hypb} Under the assumptions that $u_0\in L^1(\mathbb{H}^3)$ with mass $M>0$, and $u(x,t)$ is a solution of the heat equation in the hyperbolic space the result
\begin{equation}
\displaystyle \|u(x,t)-MP_t(x)\|_{L^1(\mathbb{H}^3)}\to 0 \quad \mbox{as \ $t\to\infty$}.
\end{equation}
is in general false, even if $u_0\ge 0$ and is compactly supported.
\end{theorem}

We give a full proof of this result  in Section \ref{sec.nonconv}. The result should be true in all dimensions $n>1$ but we prove here in dimension 3 since we make much use of explicit formulas. Explicit counterexamples are constructed. A comment on the method: we want to compare for all large times solutions with general data $u_0\in L^1(\mathbb{H}^n)$ with a multiple of the fundamental solution. We will get important information by making the comparison in the case of particular solutions that are explicitly or almost explicitly known. These solutions correspond either to a displacement of the initial mass, or to a delay in time. Both approaches give quite different results as reflected in the above theorems.

\medskip

\noindent {\sc Comments.} Following usual conventions, we use the notation $G_t(x)$ or $P_t(x)$ to denote $G(x,t)$ or $P(x,t)$ when $t$ is thought of as fixed. To avoid confusion we then use the notation $\partial_t G$, $\partial_t P$ for the partial derivative w.r.t. to time. Similar notations appear elsewhere, i.e., subscripts as variables or as derivatives, but we hope that no confusion will arise. We usually work with  $u_0\ge 0$. Since the equation is linear the results for signed data follow immediately.

\section{Preliminaries}\label{sec.prelim}

\noindent {\bf Polar geodesic coordinates.} We  recall some facts about the hyperbolic space and the heat equation posed on it. Several models are used to describe $\Hn$ in an explicit coordinate system, see \cite{BenPet92, Ratcliffe} or similar general references. Thus, one may realize $\Hn$ as an embedded hyperboloid in $\re^{n+1}$, endowed with the inherited Minkowski metric. It is also possible to use one of the two Poincar\'e realizations, so that we identify $\Hn$ with the unit ball in $\ren$ or with the upper half-space, each  endowed with an appropriate metric with the property that the Riemannian distance from any given point to points approaching the topological boundary tends to $+\infty$. Here below, it will be convenient to describe the hyperbolic space as a \sl model manifold \rm following \cite{GW79}, see also our previous work \cite{Vazhyp15}. On such a manifold, a pole $o\in \Hn$ is given, and the Riemannian metric has the form
$$
ds^2 = dr^2 + f(r)^2d\sigma^2,
$$
for an appropriate function $f$, where $r $ is the geodesic distance from the pole $o$ while
$d\sigma^2$ denotes the canonical metric on the unit sphere $\mathbb{S}^{n-1}$. The hyperbolic space $\Hn$ is obtained by making the precise choice $f(r) = \sinh r$ (hyperbolic sine function). This representation is usually referred to as  polar geodesic coordinates for hyperbolic space around the given pole. The differential expression of the Laplace-Beltrami operator $\Delta_g$ in the hyperbolic space with curvature $K=-1$ is then given by
\begin{equation}
\Delta_g u = (\sinh r)^{1-n}\frac{\partial}{\partial r}
\left((\sinh r)^{n-1}\frac{\partial u}{\partial r}\right)
+\frac{1}{(\sinh r)^2} \Delta_{\mathbb{S}^{n-1}}u\,.
\end{equation}
In such coordinates the volume element is \ $d\mu = (\sinh r)^{n-1} dr d\omega_{n-1}$, where $d\omega_{n-1}$ is the volume element on the unit sphere ${\mathbb S}^{n-1}$, see \cite{Ros98}.
 Such a formulation is most convenient in the setting of radial functions $u(r)$ or $u(r,t)$. Then, the heat equation is written as
\begin{equation}\label{he.radial}
\partial_t u =  u_{rr}''(r) + (n-1)\coth (r) u_r'(r) =
\frac1{(\sinh r)^{n-1}}\left((\sinh r)^{n-1}u_r'(r)\right)_r'\,,
\end{equation}
where $u'_r=\partial_r u$. In the sequel we will often write $\Delta$ for $\Delta_g$ if no confusion is expected.

\medskip

\noindent {\bf Fundamental solution and properties}. Let us call $P(x,o,t)$ the fundamental solution of the heat equation in hyperbolic space (h.f.s.) with center $o\in M={\Hn}$; this point  is identified with $x=0$ in the standard model representation, and also in the Poincar\'e disk. But all  points of $\Hn$ are equivalent modulo an isometrical transformation, a well-known fact. It follows that $P_t$ is only a function of the geodesic distance, $P(x,o,t)= G_t(d(o,x))$, with $d(\cdot,\cdot)$ the geodesic distance between two points. We have the representation formula for following general solutions with integrable initial data $u_0$:
\begin{equation}\label{rep.form}
u(x,t)=\int_{\Hn} u_0(y) G_t(d(x,y))\,d\mu(y).
\end{equation}
We recall that $\mu$ is the volume measure in $\Hn$ that we have just introduced. In the sequel we  perform the proofs under the further restriction that $u_0\ge 0$. For a signed solution we must only separate the positive and negative parts of the data and apply the results to both partial solutions.  We will also write $r=d(x,o)$. Some basic consequences follow from the representation formula, like conservation of mass:
$$
\int_{\Hn} u(x,t) \,d\mu(x)=\int_{\Hn} u_0(x)\,d\mu(x)\,,
$$
and the maximum principle in its strong or weak form.

\medskip

\noindent (ii) It is well-known that $G(r,t)=G_t(r)$ is a positive, smooth and decreasing function of $r>0$. Therefore, $G_t(0)$ is the maximum value of $G_t(r)$, $r\ge 0$. In this way, we obtain a $L^1$-$L^\infty$ smoothing effect
\begin{equation}\label{upper.fs}
\|u(\cdot,t)\|_\infty\le G_t(0)\|u_0\|_1.
\end{equation}
We see that the fundamental solution gives the worst case in this estimate.

\medskip

\noindent (iii) The fundamental solution of the heat equation in hyperbolic space has a clear explicit form in dimension $n=3$. It reads:
 \begin{equation}
 G(r,t)= C_3t^{-3/2}e^{-t}\frac{r}{\Sh r} e^{-r^2/4t}, \quad C_3=(4\pi)^{-3/2}.
 \end{equation}
 In other dimensions $n\ge 2$ the formulas are rather involved, see \cite{Grig98}.  There is a very useful recurrence relation that allows to derive the fundamental solution for all dimensions from just two of them.
 It is contained in \cite{DavMand}.

  \begin{proposition}\label{prop.recur} If $P_{(n)}(x,t)$ is the fundamental solution around a given pole in dimension $n\ge 1$ expressed in polar geodesic coordinates as $G_{(n)}(r,t)$, then
\begin{equation}
  v(r,t)=-e^{-nt}\,(\mbox{\rm \Sh} r)^{-1}\partial_r G_{(n)}(r,t)
\end{equation}
  gives the expression of the fundamental solution around the same pole in dimension $n'=n+2$, but for a constant. We have \ $v(r,t)= 2\pi \,G_{(n+2)}(r,t)$.
  \end{proposition}

In this way, formulas for $G_{(n)}(x,t)$ can be derived for all odd $n$ from $G_{(3)}(x,t)$. Very precise information about function $G_t$ in all dimensions is given by the following  estimate by Davies \cite{Davies}, see also \cite{DavMand}.

\begin{proposition}\label{prop.davies}   For all $n\ge 1$ there exists a positive constant $c_n$ such that
\begin{equation}\label{est.h}
\frac1{c_n} h_n(r, t) \le G_{(n)}(r, t) \le c_n h_n(r, t)
\end{equation}
for all $t > 0$ and $r > 0, $ where
\begin{equation}\label{est.h2}
h_n(r, t) =\frac1{(4\pi t)^{n/2}}\,e^{-\lambda_1 t-\frac{n-1}{2} r- \frac{r^2}{4t}}
 (1 + r + t)^{(n-3)/2}(1 + r).
\end{equation}
and $\lambda_1=(n-1)^2/4.$
\end{proposition}

Note that $\lambda_1$ is the bottom of the spectrum of the Laplacian in $\mathbb{H}^n.$ It is often useful to write the exponential term as
$$
\displaystyle e^{-\displaystyle\lambda_1 t-\displaystyle\frac{n-1}{2} r- \displaystyle\frac{r^2}{4t}}= \ e^{-  \displaystyle\frac{\left(r+(n-1)t\right)^2}{4t}}.
$$

If we fix $r= 0$ then we have   \ $ G(0, t)\approx t^{-n/2}$ \
 as $t \to 0$. Actually, the fundamental solution of hyperbolic space looks like the Gaussian kernel of Euclidean space as $t\sim 0$, but we will not be much interested in small times. On the other hand, in the limit $t\to \infty$ we have a very different behaviour:
\begin{equation}
G(0, t)\approx t^{-3/2}e^{-\lambda_1 t}.
\end{equation}
In this way we obtain a precise $L^1$-$L^\infty$ smoothing effect.

\begin{proposition}\label{prop.sm.eff} There is a constant $C_n>0$ such that for all integrable solutions $u(x,t)$ the following decay estimate is valid for all $t\ge 1$
\begin{equation}\label{reg.eff}
\|u(\cdot,t)\|_\infty\le C_n t^{-3/2}e^{-\lambda_1 t}\|u_0\|_1\,.
\end{equation}
The exponents in the right-hand expression as $t\to \infty$ are optimal. This is valid for all $n>1$.\end{proposition}

\section{Pointwise bounds and Harnack-like inequalities}
\label{sec.pbhi}

Now that we have some precise knowledge  about the behaviour of the fundamental solution, we can use it as a comparison tool for general solutions with nonnegative, compactly supported data thanks to the following result.

\begin{theorem}\label{thm.pb} For every solution $u(x,t)$ with initial data $u_0\ge 0$ having mass $M>0$ and supported in a ball of radius $R>0$ around the pole $o$, and for every $L>0$ there is a constant $C=C(R,L)>1$ such that the following comparison holds
\begin{equation}\label{pbwrtfs}
\frac1{C} M\,P_t(x)  \le u(x,t)\le C M\,P_t(x)
\end{equation}
for all $r=d(x,o)\le Lt$ and $t\ge 1$.
\end{theorem}

\noindent {\sl Proof.} By linearity we may assume that $M=1$.

\noindent (i) We start with the upper bound. The proof of the bound will be based on the fact that there is a worst case  for all solutions with $u_0$ supported in the ball of radius $R$ since when take a point $x$ located far away from $o$ (more distance than $R$) in the direction of the first axis, $x=(r,0,\dots,0)=r\,e_1$, inspection of the representation formula \eqref {rep.form} shows that we can move all the mass to the point $y_1$ with  coordinates $x_R=(R,0,\dots,0)$ and then since $d(x,x_R)=d(x,0)-R<d(o,x)=r $ and the fact that $G_t$ is decreasing in $r$ we get for $x=r\,e_1$, $r>R$,
$$
u(x,t)<u_1(x,t)= G_t(d(x_R,x))=G_t(r-R).
$$
Since the equation is invariant under rotations, considering points of the form $x=r\,e_1$, $r>0$, is no restriction. This allows us to obtain  desired upper bound in the annular region
 $$
 \mathcal{A}=\{(x,t):   \, R<d(x,o)=r<Lt, \, t>0\}.
 $$
Together with the inequality  $ u(x,t)\le G(r-R,t)$, we use the estimates \eqref{est.h}, so that
 $$
 G(r-R,t)\le c_1(n)\,h_n(r-R,t)\,.
 $$
Using the explicit formulas for $h_n$ we observe that
 $$
 \frac{h_n(r-R,t)}{h_n(r,t)}\le e^{(n-1)R}e^{rR/2t}\left(1-\frac{R}{1+r}\right)\left(1-\frac{R}{1+r+t}\right)^{(n-3)/2}.
 $$
 Therefore, in the region $\mathcal{A}$ we get
 $$
 u(x,t)\le c_1(n)c_2(R,L)\, h_n(r,t)\le c_3(n,L,R) \, G(r,t).
  $$

\noindent  (i') We still need to estimate the solution  inside the ball of radius $R$. But we know from  \eqref{upper.fs} that for every $(x,t)$
 $$
 u(x,t)\le  P(0,t).
$$
Changing $P(0,t)$ into $0<P(r,t)$ with $r<R$ and $t>1$ implies inserting a constant depending on $R$. Putting both estimates together, the upper bound in \eqref{pbwrtfs} follows.

\noindent (ii) We remark that under the assumptions and notations of (i) a similar argument applies to the first step of the lower bound. Now the worst case then consists of moving the initial mass to $x_{-R}=(-R,0,\dots,0)$ and then we get
$$
u(x,t)>u_2(x,t)= G_t(d(x_{-R},x))=G_t(r+R)\ge c_n^{-1}h_n(r+R,t).
$$
for all $x=r\,e_1$, $r>0$. ,
 $$
 \frac{h_n(r+R,t)}{h_n(r,t)}\ge e^{-(n-1)R}e^{-R^2/4t}e^{-rR/2t}\left(1+\frac{R}{1+r}\right)\left(1+\frac{R}{1+r+t}\right)^{(n-3)/2}.
 $$
The last quantity is bounded below for all $0<r<Lt$ and $t>1$. \qed

 \medskip

\noindent \bf Remarks. \rm 1) The control of the tail $r\ge Lt$ by a similar expression is  false, but it is true if we insert time delays into $P$.

2) Note that in the Euclidean space setting this result is true  with $C\to 1$ as $t\ge t_0\to\infty$,  but such a fine convergence is false in the hyperbolic space, as we will show below.

\medskip

We easily derive from the above result a Harnack inequality of  elliptic type for such a class of solutions.

\begin{corollary}\label{cor.hi} For every two solutions $u_i(x,t)$, $i=1,2$, with initial data $u_{0i}\ge 0$ having mass 1 and supported in a ball of radius $R>0$ around the pole $o$, and for every $L>0$ there is a constant $C=C(R,L)>1$ such that the following comparison holds
\begin{equation}\label{pbwrtfs}
\frac1{C} u_1(x_1,t)  \le u_2(x_2,t)\le C u_1(x_1,t)
\end{equation}
for all $r_i=d(x_i,o) \le Lt$ and $t\ge 1$.
\end{corollary}

\section{Mass analysis for the fundamental solution}
\label{sec.mass}

We now proceed with the study of the long-time behaviour of solutions. We return in this section to the fundamental solution. As an important tool we introduce the concept of {\bf mass function}. It is defined for any integrable radial distribution $f(r)\ge 0$ in hyperbolic space by the expression
\begin{equation}\label{mass.fn}
M(f)(r)=\omega_n\int_0^r (\Sh r)^{n-1}f(r)\,dr\,,
\end{equation}
where $\omega_n$ is the measure of the unit sphere $\mathbb{S}^{n-1}$.
Without loss of generality we may assume that the total mass equals 1 so that $M(f)(r)$ will be a monotone function of $r$ with $M(f)(0)=0$, $M(f)(\infty)=1$.

If we apply this concept to an evolving distribution $u(r,t) $ that satisfies the heat equation, we obtain an evolving mass function $M(r,t)$  such that $\partial_r M(r,t)=\omega_n(\Sh r)^{n-1}u(r,t)$. From equation \eqref{he.radial} we derive the evolution equation for the mass function $M(r,t)$:
\begin{equation}\label{mass.eqn}
\partial_t M=\omega_n(\Sh r)^{n-1}\partial_r\left(\partial_r M/(\Sh r)^{n-1}\right)= M_{rr} - a(r) M_r,
\end{equation}
which is a modified 1D heat equation with a drift term that has variable speed:
$$
 a(r)=-(\Sh r)^{n-1}((\Sh r)^{1-n})_r=(n-1)\,\Cth r\,,
 $$
moving in the forward direction. In the sequel we will use the fact that $a(r)\to (n-1)$ as $t\to\infty$.  We may also find a similar equation for the radial mass density $$\rho(r,t)=\partial_r M(r,t)=\omega_n(\Sh r)^{n-1}u(r,t),
$$
see details below. We recall that here and in the sequel $r=d(x,o)\ge 0$.


\subsection{Location of the mass bulk. Linear rate of mass escape for $n=3$}
When studying the heat equation in  hyperbolic space it is important  to find out where the most of the mass is actually located. We recall that in Euclidean space the mass of the fundamental solution spreads with time at distances of the order of $O(t^{1/2})$. This is the so-called Brownian scaling. There is a big difference in hyperbolic space due to the presence of a drift term.  We prove the following result about the mass location of the fundamental solution for the heat flow in hyperbolic space that may look very striking at first sight. We work in $n=3$ for the moment in order to exploit the extra sharpness given by the explicit formulas.

 \begin{theorem} Let $n=3$. For all large $t\gg 1$ most of the mass of the  fundamental solution is located in a thin annulus of the form $\{ r_1(t)\le r=d(x,o) \le r_2(t)\}$, where
\begin{equation}
 r_1(t)=2t-kt^{1/2}, \quad  r_2(t)=2t+kt^{1/2}\,,
\end{equation}
 with $k>0$ large but constant in time.
 \end{theorem}

 \noindent {\sl Proof.} The mass of the fundamental solution in a thin annular region of the form $\{ r_0<r<r_0+dr\}$ is equivalent to
$$
dM(r_0)=C (\Sh r)^2 G_t(r)dr =
C\,\,t^{-3/2} r \frac{\Sh r}{e^r} e^{-\frac{(r-2t)^2}{4t}}dr.
$$
This is better understood in terms of the variable $\xi=(r-2t)/t^{1/2}$. We get for large $t$
$$
dM=C\,\,\frac{r}{t}\frac{\Sh r}{e^r} e^{-\xi^2/4}d\xi.
$$
After a quick examination we conclude that for large $t$  there exists a significant part of the total mass inside the region where $\xi$ is bounded, $|\xi|\le k$, which implies that $r\sim 2t+ O(t^{1/2})$ (remember that $r=2t+ \xi t^{1/2}$). Note that
$$
 \frac{r \Sh r}{te^r} \to 1
  $$
 as $t\to\infty$ with $r/t\to 2$. In particular, the mass located outside the interval $I_k=\{-k\le \xi\le k\}$ is very small, of the order of the error function.   \qed

Note that the so-called \sl mass line \rm $r=2t$ refers actually to a expanding spherical surface when seen as a set in $\mathbb{H}^3$, or a conical hypersurface when seen as a subset of space-time, $\mathbb{H}^3 \times (0,\infty)$.

\medskip

 \noindent {\bf Other geometrical lines.}
 We want to relate this escape line with other geometrical or mechanical distances. We may define {\sl the half-mass  line} as the line $r=r_m(t)$  so that
 $$
 \int_0^{r_m(t)} dM(r)=1/2\,.
 $$
Our formulas above show that $r_m(t)/2t\to 1$ as $t\to \infty$, more precisely
 $r_m(t)=2t + o(t^{1/2})$. Thus, for large $t$ the mass is mostly located in a $t^{1/2}$-neighbourhood of the  half-mass  line.

 A related distance is given by the {\sl sign-change line}, $r= r_{s}(t)$, where $\partial_t G(r,t)=0$ (with $\partial_t G(r,t)<0$ for $0\le r<r_{s}(t)$). An easy computation  gives the exact value
 $$
 r_{s}^2 (t)= 6t+4t^2.
 $$
We have $r_s(t)=2t+ O(t^{1/2})$ as $t\to \infty$. Note that the approximation $r_s^2(t)\approx 6t$ is good for small times and exactly represents the Euclidean situation.

\medskip

\subsection{Escape analysis for general dimensions. Ballistic motion} The same analysis in hyperbolic space with $n>1$, $n\ne 3$ leads to the mass line \ $r=(n-1)t$; the mass is concentrated for large $t$ around that line  in a region with local width $O(t^{1/2})$.  We prove this by using the upper estimates in Proposition \ref{prop.davies}. Indeed, we can prove that the tail integral
 $$
I_{+}(L;T) = \int_{r(L,t)}^\infty t^{-n/2}\,e^{-(r-(n-1)t)^2/4t} (1 + r + t)^{(n-3)/2}(1 + r)\,dr
$$
with $r(L,t)=(n-1)t+ L\,t^{1/2}$, is small for large fixed $L$ and all large $t\gg 1$. But for large $t$ and putting $r=(n-1)t+\xi t^{1/2}$ we get:
$$
I_{+}(L;T) \le  C\int_{L}^\infty t^{-(n-1)/2}\,e^{-\xi^2/4}(t + \xi t^{1/2})^{(n-1)/2}\,d\xi\le
\int_{L}^\infty e^{-\xi^2/4}(1 + \xi t^{-1/2})^{(n-1)/2}\,d\xi
$$
which goes to zero with $L\to \infty$, uniformly  in $t$. The same applies to the near field remainder
$$
I_{-}(L;T) = \int_0^{r_{-}(L,t)} t^{-n/2}\,e^{-(r-(n-1)t)^2/4t} (1 + r + t)^{(n-3)/2}(1 + r)\,dr
$$
with $r_{-}(L,t)=(n-1)t- L\,t^{1/2}$.

\medskip

\noindent {\bf Note.} We have learnt the recent work of Lemm and Markovic  who prove in Appendix A of \cite{LemMark2018} the drifted location of the heat kernel that we have described. The feature is described by them as {\sl ballistic propagation}.

\subsection{Dimensional explanation} The escape of the mass with a precise linear motion looks surprising but there are simple dimensional considerations to support it (in rough terms). First of all, the well-known fact that $\lambda_1=(n-1)^2/4>0$ is the first eigenvalue of the Laplace-Beltrami operator in $\Hn$ explains the appearance of the exponentially decreasing factor $ e^{-\lambda_1 t} $ in the sup norm size of the h.f.s. The factor $t^{-3/2}$ serves as a minor correction and is related the Euclidean case.

Next, since there is conservation of mass too,  the h.f.s needs to spread at least in a region of volume inverse to the sup norm, i.e. roughly $V\sim e^{\lambda_1 t}$ (we forget the power correction in rough approximation). Since the volume of a ball grows like  $V(B_r(0))\sim (\Sh r)^{n-1}$, we find after this  rough calculation that most of the mass is located at or beyond \
 $r_1\sim {\lambda_1} t/(n-1)=(n-1)t/4$, which is our result up to a factor $1/4$.
 Of course, this does not explain why the mass concentrates mostly near the boundary of a precise ball.


\section{The 1D asymptotic limit. 1D heat equation with drift}\label{sec.1Ddrift}

  We have established in Section \ref{sec.mass} the asymptotically linear {\sl escape motion} of the mass bulk of the h.f.s.: for large $t$ the mass is mostly located next to the line $r=(n-1)t$. This phenomenon clearly does not happen in the Euclidean case. Next, we want to see the detailed evolution near this line. Actually,  when $t\to \infty$, we will show that the fundamental solution becomes a  standard 1D heat equation with constant drift (after proper renormalization).

 Let us prove this fact in $n=3$ by direct inspection: if $\rho(r,t)=\partial_r M(r,t)$ is the mass density of the h.f.s., then in three space dimensions we get
$$
\rho(r,t)= 4\pi (\Sh r)^2\, G(r,t)=(4\pi t)^{-1/2}\frac{r\, \Sh r}{t e^r} e^{-(r-2t)^2/4t}.
 $$
 As $t\to\infty$ and for $r/2t\to 1$ we get the equivalence
 $$
\rho(r,t)\sim E_1(s,t)= (4\pi t)^{-1/2} e^{-s^2/4t}\,,
$$
where $E_1(s,t)$ is the fundamental solution of the heat equation in the real line $\re$, and $s=r-2t$ is the radial coordinate w.r.t the mobile frame. We  easily get the following precise result.

\begin{theorem}\label{conv.fsto1d} If $G(r,t)$ is the fundamental solution of the heat equation in hyperbolic space $\mathbb{H}^3$ and $E_1(r,t)$ is the fundamental solution of the heat equation in the real line $\re$, we get the limit \rm
\begin{equation}
\lim_{t\to\infty}
\frac{4\pi\, (\Sh r)^2\, G(r,t)}{E_1(r-2t,t)}=1
\end{equation}
\sl uniformly on sets of the form:  $|s|=|r-2t| \le L t^{1/2}$, $L>0$.  Moreover, we have the weighted $L^1$-error estimate \rm
\begin{equation}
\lim_{t\to\infty} \int_0^\infty |4\pi(\Sh r)^2\, G(r,t)-E_1(r-2t,t)|\,dr=0\,.
\end{equation}
\end{theorem}

\subsection{Limit behaviour in general dimensions}

This convergence result can be generalized to general space dimensions $n\ge 2$ . We will do it via the PDE satisfied by $\rho$, since using the same approach with direct computations soon becomes cumbersome. Thus, from \eqref{mass.eqn} we get the following equation for the mass density
$\rho(r,t)=\partial_r M(r,t)=\omega_n(\Sh r)^{n-1}G_t(r)$:
\begin{equation}\label{massdens.eqn}
\partial_t\rho +(n-1)\partial_r\rho= \partial^2_{rr}\rho - (n-1)\partial_r(b(r)\rho)\,, \quad
b(r)= \Cth(r)-1\,.
\end{equation}
As  $t\to\infty$ and for $r/(n-1)t\to 1$, we see that $b(r)\to 0$ (very quickly) so the last term is really lower order. Hence, we hope to get convergence of $\rho$ to a solution of the limit equation:
 \begin{equation}\label{massdens.eqnlim}
\partial_t\rho +(n-1)\partial_r\rho= \partial^2_{rr}\rho\,,
\end{equation}
which is a heat equation with the expected outgoing drift. The limit solution $\overline{\rho}$ must be a 1D Gaussian for large times after eliminating the drift. Complete details are needed to justify the proof. They go as follows:

\noindent (i) We introduce new coordinates to correct for the drift term
$$
 s=r-(n-1)t, \quad \widehat{M}(s,t)={M}(r,t), \quad \widehat{\rho}(s,t)=\rho(r,t)\,.
$$
Then the equations for mass and density become resp.
$$
\partial_t \widehat{M}= \partial^2_{ss}\widehat{M} -\widehat{b}(s) \,\partial_s\widehat{M},
$$
with $\widehat{b}(s)= \Cth(r)-1>0$, and
$$
\partial_t\widehat{\rho} = \partial^2_{ss}\widehat{\rho} - (n-1)\,\partial_s(\widehat{b}(s)\widehat{\rho})\,.
$$
Both equations are posed in the expanding parabolic domain
$$
Q=\{(s,t): \ t>0, \ -(n-1)t < s< \infty\}.
$$

\noindent (ii) We remark that the last term in both equations can be considered as an asymptotically small perturbation with respect to the heat equations satisfied approximately by $\widehat{M}$ and $\widehat{\rho}$. This will imply convergence of both as $t\to\infty$ to Euclidean limit for all $s=r-(n-1)t$. The asymptotic analysis of evolution processes that can be viewed as asymptotically small perturbations of known dynamical systems has a long tradition in PDEs and Mechanics. We refer to the book \cite{GV2004} for the theory and the application to a number of nonlinear heat equations. The main result says that, under rather simple given assumptions, the orbits of the perturbed system converge as $t\to\infty$ to the $\omega$-limit (i.e., the asymptotic states) of the master system. Here, the master equation is the 1D heat equation in the line, and its $\omega$-limit in the class of integrable solutions is formed by multiples of the Gaussian kernel, a well established fact, \cite{VazGN}.

\noindent (iii) The convergence proof in the present case may proceed as follows.  In view of the fact that the tails are under uniform control, we need interior estimates, specially in the range where $r\sim (n-1)t$. First of all,
$$
\partial_s \widehat{M}= \partial_r  M(r,t)=\rho(r,t)=c(\Sh r)^{n-1}u(r,t)\,,
$$
which is positive and bounded. According to the estimate of Proposition \ref{prop.davies}, it is bounded uniformly for large $t$, and of order $t^{-1/2}$ for $r\sim t$.
Furthermore, $M_t\le 0$ and
$$
-M_t=-\int_0^r (\Sh r)^{n-1}u_t\,dr=-\int_0^r ((\Sh r)^{n-1}u_r)_r\,dr=
-(\Sh r)^{n-1}u_r\ge 0,
$$
which for $r\sim (n-1)t$ gives
$$
-M_t= e^{nt}(\Sh r)^{n}u_{(n+2)}(r,t)\le C e^{nt}(\Sh r)^{-1}re^{-(2t)^2/4t}t^{-3/2}\sim Ct^{-1/2}.
$$
Now, $\partial_t\widehat{M} =M_t - (n-1)M_s$, so it also has a good bound as $t\to\infty$.

\noindent (iv) We want to pass to the limit $t\to\infty$ in a typical heat equation way. For that we perform a rescaling on the space variable and  consider the convergence of  the mass expressed in the rescaled variable, $N(\xi,t)= \widehat{M}(\xi t^{1/2},t)$, i.e., as a function of $\xi=s/t^{1/2}$ and $t$.  Using compactness and the far field estimates, we conclude that \ $N(\xi,t+t_k)$ converges along subsequences $t_k\to\infty$  to a solution of the rescaled heat equation, $\overline{N}$, the limit is monotone in $\xi$ and joins the level $\overline{N}(-\infty,t)=0$ to $\overline{N}(-\infty,t)=1$. A close analysis of the limit as solution of the rescaled heat equation shows that $\overline{N}$ must be the integral of the Gaussian kernel (i.e., its accumulated mass). In other words,
\begin{equation}\label{mconv}
\lim_{t\to\infty} |\widehat{M}(s,t)-M_1(s)|\to 0\,,
\end{equation}
uniformly in $s\in \re$. $M_1$ denotes the 1D mass function of the stationary Gaussian kernel.

\noindent (v) We continue with the convergence of
$$\overline{\rho}=\partial_\xi N(\xi,t)=t^{1/2}\partial_s \widehat{M}=t^{1/2}\widehat{\rho}.
$$
 A previous estimate shows that $\overline{\rho}$ is a bounded function of $\xi$ and $t$ for large $t$. Also, $\overline{\rho}_\xi=t^{1/2}\widehat{\rho}_s(s,t)$, hence
$$
\overline{\rho}_\xi =t^{1/2}{\rho}_r=Ct^{1/2}((\Sh r)^{n-1}u)_r=Ct^{1/2}(\Sh r)^{n-1}u_r + C(n-1)t^{1/2}(\Sh r)^{n-2}(\Ch r) u.
$$
We conclude from previous estimates that $\overline{\rho}(\xi,t)$ is compact and therefore we get the convergence in $L^1$ to the derivative of $M_1(\xi)$, i.e., the Gaussian kernel, by virtue of \eqref{mconv}. In this way we get the generalization of Theorem \ref{conv.fsto1d} to  $n$ dimensions.

\begin{theorem}\label{conv.fsnto1d} If $G(r,t)$ is the fundamental solution of the heat equation in hyperbolic space $\mathbb{H}^n$, $n\ge 2$, and $E_1(r,t)$ is the fundamental solution of the heat equation in the real line $\re$, we get the limit
\begin{equation}\label{rhoconv}
\lim_{t\to\infty} \int_{\xi_0}^\infty |\overline{\rho}(\xi,t)-E_1(\xi,1)|\,d\xi=0\,,
\end{equation}
 where, $\xi_0=-(n-1)t^{1/2}$ means $r=0$. In other words,\rm
\begin{equation}
\lim_{t\to\infty} \int_0^\infty |\omega_n(\Sh r)^{n-1}\, G(r,t)-E_1(r-(n-1)t,t)|\,dr=0\,.
\end{equation}
\end{theorem}

\section{Long-time convergence results for radial solutions}
\label{sec.rad}

We want to prove the following main result.

\begin{theorem} \label{thm.radconv} Let $u(r,t)$ a positive solution of the heat equation in $\Hn$ with radially symmetric data having compact support and mass 1, and let $P(r,t)$ be the fundamental solution. Under these assumptions we get
\begin{equation}\label{rd.conv1}
\lim_{t\to\infty}\|P(r,t)-u(r,t)\|_{L^1(\mathbb{H}^n)}=0,
\end{equation}
\end{theorem}

The condition of compact support on the data may be replaced by  strong decay at infinity, but we will not discuss that issue.

\subsection{Convergence with rate for a solution with time delay for $n=3$}

We start the proof  dimension 3 since the details are quite precise. We compare  the fundamental solution $P(r,t)$ with its time-delayed version $Q(r,t)= P(r,t+1)$ for large $t$, as a stepping stone for comparison of more general radial solutions.

\noindent (i) Preliminary calculation. We have
$$
\log P= \log C- (3/2) \log t+  \log (r/\Sh r) -t -\frac{r^2}{4t},
$$
so that
$$
\frac{\partial_t P}{P}= -\frac3{2t}-1+\frac{r^2}{4t^2}.
$$
We immediately see that
$$
\frac{\partial_t P}{P}\ge  -\frac3{2t}-1.
$$
This proves that for large $t$ we have
$$
P(r,t+1)\ge P(r,t)e^{-1}\left(t/(t+1)\right)^{3/2}.
$$
which relates the mass of $P_t$ and $P_{t+1}$ but for a factor $e$. This is not enough for our purposes.

To continue we put  $r=2t+ \xi t^{1/2}$. Then $r^2=4t^2+4\xi t^{3/2} + \xi^2 t$, and
$$
\frac{\partial_t P}{P}= -\frac3{2t}-1+\frac{r^2}{4t^2}=-\frac3{2t} + \xi t^{-1/2}+ \frac14\xi^2 t^{-1},
$$
so that we have $|\partial_t P|\le C P\, t^{-1/2}$ for bounded $\xi$: $|\xi|\in L $, uniformly in $t\gg 1$.

\medskip

\noindent (ii) Compare now the solutions $P(r,t)$ with $Q(r,t)= P(r,t+1)$ for large $t$. In view of our previous analysis of the analysis of the location of the mass, we can use (i) to prove that the local mass of the difference $|P_t-Q_t|$ is very small for large $t$
 with an order of convergence. We know that the zone of large mass of solution $P$ is $\{-L<\xi<L \}$ for some large $L$, hence $\{ 2t -Lt^{1/2}<r<2t +Lt^{1/2}\}$. The same is true for $Q$ in the region
$$
2(t+1)-L(t+1)^{1/2} <r< 2(t+1) +L(t+1)^{1/2}=2t+ Lt^{1/2} +2 + O(t^{-1/2}).
$$
Both regions are very similar for large $t$. In the intersection of both regions we have a very small relative variation of $P$, hence a very small variation of the mass with unit time. In the complement of this region, we get a very small mass for both solutions anyway. Putting these facts together, we get the $L^1$ convergence result.

\noindent {\bf Note.} We can get a decay rate by estimating the mass outside in terms of the bound $L$ for $\xi$. This is to be done, but not difficult.


\subsection{Convergence example in general dimensions}

We can use the asymptotic result of Theorem \ref{conv.fsnto1d} to prove that the $L^1$ convergence of $P_t(x)-P_{t+1}(x)$ to zero in $L^1(\mathbb{H}^n)$ as $t\to\infty$. Indeed, from \eqref{rhoconv} we get
\begin{equation}\label{rhoconv.del}
\lim_{t\to\infty} \int_{-\infty}^\infty |\overline{\rho}(\xi,t+1)-\overline{\rho}(\xi,t)|\,d\xi=0\,,
\end{equation}
where the functions are extended by zero for $r<0$ for ease of notation. This means that
\begin{equation}
\lim_{t\to\infty} \int_{\Hn} |P(x,t+1)-P(x,t)|\,d\mu(x)=0\,.
\end{equation}

No rate is given in this case.


\subsection{General radial data}

The same positive conclusion holds for a wide class of nonnegative radial initial data. We use intersection number and mass comparison. In order to do the analysis we use polar geodesic coordinates for functions $u(r,t)\ge 0$.

(i) We recall the mass function \eqref{mass.fn} and the mass equation \eqref{mass.eqn}.
For smooth radial solutions of the hyperbolic heat equation, the mass equation satisfies the maximum principle in the sense that  when $u_1$ and $u_2$ are smooth solutions of the heat equation with total mass 1, and we have the comparison $M_1(r,0)\le M_2(r,0)$ for all $r>0$ for the respective mass functions, then we have
$$
M_1(r,t)\le M_2(r,t)\quad \mbox{for all \ } r,t>0.
$$

(ii) We now recall the basic idea of intersection comparison, also called lap number. The number of intersections of two solutions $I(u_1,u_2) $ of the heat equation does not increase with time. In particular, if at $t=0$ the number is 1 then for later times it must be 0 or 1. But if both have mass 1, it cannot be zero intersections. Hence, it is 1 all the time.

\noindent {\bf Note.} The technique works in one dimension, or for radially symmetric solutions in several dimensions, and is based on counting the evolution in time of the ``number of intersections of two solutions'', a rough idea that can be made precise with the names  {\sl intersection number} or {\sl lap number}. These concepts have been investigated in works by Sattinger \cite{Sat69}, Matano \cite{Mat82},  Angenent \cite{An88} and others, and were used by Galaktionov and the author in a number of papers, cf. \cite{GV2004}.

(iii) We  apply both principles to a solution $u$ that is sandwiched between $u_1(r,t)=P(r,t)$, the fundamental solution, and $u_2=P(r,t+1)$. The precise conditions are
$$
M_1(r,0)=1\ge M(r,t)\ge M_2(r,t).
$$
and
$$
I(u_1,u)(t=\ve)=1, \quad I(u,u_2)(t=0)=1,
$$
with the intersection such that $u-u_1$ and $u_2-u$ are positive for all large $r$.
Note that the last condition also happens for $I(u_1,u_2)$.
These conditions will be conserved for all times $t>0$. Let us fix $t$ large so that
by the preceding theorem $M_1(r,t)-M_2(r,t)\le \ve$ for all $r$. Let $r_* $ be intersection point of $u_1$ and $u$. We have
$$
\int (u_1-u)_+\,d\mu=C(M_1(r_* )-M(r_* ))\le C(M_1(r_* )-M_2(r_* )) \le C\ve,
$$
while conservation of mass implies that
$$
\int (u-u_1)_+\,d\mu=- \int (u_1-u)_+\,d\mu.
$$
This ends the proof that Theorem \ref{thm.radconv} holds for $P$ and $u$.

(iv) for more general radial integrable data we use approximation. \qed

\subsection{Sup and $L^p$ convergence} 

We state now the convergence in the $L^p$ norms. The convergence in the sup norm is not so strong as the following result shows:

\begin{corollary} Under the assumptions of Theorem \ref{thm.radconv} we have
\begin{equation}\label{rd.conv2}
\|P(r,t)-u(r,t)\|_{L^{\infty}(\mathbb{H}^n)}=O(t^{-3/2}e^{-\lambda_1 t})\quad
\mbox{as $t\to\infty$}\,,
\end{equation}
as well as the optimal result
 \begin{equation}\label{rd.convp}
\|P(r,t)-u(r,t)\|_{L^{p}(\mathbb{H}^n)}=o(t^{-3(p-1)/2p}e^{-\lambda_1(p-1) t/p})\quad
\mbox{as $t\to\infty$}.
\end{equation}
valid for all \ $1<p<\infty$,
\end{corollary}

\noindent {\sl Proof.} (i) The sup result is a consequence of Theorem \ref{thm.pb}. It also follows from the convergence \eqref{rd.conv1} and the regularizing effect \eqref{reg.eff}. Indeed, given $\varepsilon>0$ for some large $T$ we have
$$
\|P(r,T)-u(r,T)\|_{L^1(\mathbb{H}^n)}\le \varepsilon,
$$
so that for $t>T$ we get
$$
\|P(r,t)-u(r,t)\|_{L^\infty(\mathbb{H}^n)}\le \varepsilon e^{\lambda_1 T}t^{-3/2}e^{-\lambda_1 t}= Ct^{-3/2}e^{-\lambda_1 t}.
$$
(ii) The result for $1<p<\infty$ follows by interpolation.

\medskip

\noindent {\bf Remark.} The  result for $p=\infty$ is weaker. There is a reason for that.
If we take the explicit h.f.s in $n=3$ and perform a delay $u(x,t)=G_{t+a}(r)$ we see that
$$
\frac{u(0,t)}{G_t(0)}= e^{-\lambda_1 a}\frac{t^{3/2}}{(t+a)^{3/2}}\to e^{-\lambda_1 a}\,,
$$
so that
$$
t^{3/2}e^{\lambda_1 t} \|P(r,t)-u(r,t)\|_{L^{\infty}(\mathbb{H}^n)}
$$
does not to go zero as $t\to\infty$.


\section{Horospheric solutions}\label{sec.horo}

There is a class of solutions of the heat equation in $\Hn$  that also enjoy a special symmetry, called horospheric symmetry, that allows for an easy analysis. They are interesting examples  and serve to derive consequences for the general theory. In particular, they allow to demonstrate the phenomenon of mass escape, even the precise speed $n-1$, and actually it does  in a very clear way. The class is best presented  in the upper half-space representation of the hyperbolic space, identifying  $\Hn$  as $\re^{n-1}\times \re_+$ (one of the mentioned Poincar\'e representations). We take as coordinates for  points in $\Hn$  the pairs $(x,y)$ with $x\in\ren$ and $y>0$, and then the metric is given by $ds^2=y^{-2}(dx^2+dy^2)$.

We consider only solutions of the heat equation that depend only on the variable $y$. Since the lines $y=c$ are so-called horospheres (horocycles in 2D), we call these solutions $u(y,t)$ {\sl horospheric solutions}. The heat equation is then written as a modification of the 1D equation
\begin{equation}
u_t=y^2 u_{yy}-(n-2)y\,u_y.
\end{equation}
There is a very handy transformation to transform it into the standard 1D heat equation. We just put $u(y,t)=v(z,t)$, with $ z=\log(y)$, and get
\begin{equation}
v_t=v_{zz}-(n-1)v_z.
\end{equation}
This is the same equation with drift that we have found before as limit of the  equations that govern the class radial solutions for $t\gg 1$, when the geometric factor, $\Cth r$, tends to just 1, cf.  equation \eqref{massdens.eqnlim}. Note also that $z=\int dy/y=\int ds$ is just the geodesic distance when going along the $y$ axis, normal to the horospheres.

Hence,  if we take horospheric solutions with  initial data $u_0(y)\ge 0$  and satisfying the 1D integral condition:
\begin{equation}\label{1D.int}
\int_0^\infty u_0(y)y^{-1}dy=\int_{-\infty}^\infty u_0(e^z)dz<\infty,
\end{equation}
the asympotics as $t\to\infty$  is given by an expanding Gaussian $E_1(s,t)$ with respect to the moving space variable $s=z-(n-1)t$, and the bulk moves sideways in the forward direction as function $v(z,t)$, with constant speed $(n-1)$ measured in the $z$ scale.

\begin{theorem} With the previous notations and assumptions, we have
 \begin{equation}
\lim_{t\to\infty} \left|t^{1/2}v(z,t)-E_1((z-(n-1)t)/t^{1/2})\right|=0\,
 \end{equation}
uniformly in $s$.
\end{theorem}

The reader is asked to picture the equivalent drift motion expressed with respect to the original $y$ variable.

\noindent {\bf Comments.} 1) In the way of comparison of the class of radial solutions and horospheric solutions, let us recall that the integrability condition \eqref{1D.int} we have imposed on the latter does not imply integrability in $\Hn$ since we are not integrating in the variables $(x_1,...,x_{n-1})$. Therefore, the $\Hn$ integral is infinite.

\noindent 2) A consequence of this different integrability is that our horospheric solutions decay like
$$
u(y,t)=v(s,t)=O(t^{-1/2})\,,
$$
as time goes to infinity. The exponential time factor of the radial solutions is canceled. On the other hand, the drift has the same size, only dependent of the dimension.

\section{Negative asymptotic convergence results for displaced masses}
\label{sec.nonconv}

This study has two parts, one about the relative pointwise error that is done near the origin, and implies the sup norm results of Theorem \ref{thm.conv.hypb}. Another part deals with the examples exhibiting  the lack of convergence of the $L^1$ error, that is shown to happen at the critical distance $d(x,o)\sim (n-1)t$.

\subsection{The pointwise error }

\noindent (i) Let $P_t(x,o)=G_t(r)$, the fundamental solution centered at $x=o$. We want to prove that
$$
\lim_{t\to\infty}t^{3/2}e^{\lambda_1 t} P_t(x;o)= Q_n(r)>0
$$
for a certain $C^\infty$ function $Q_n$ that is bounded and monotone nonincreasing  in $r$. From Proposition \ref{prop.davies} we know that
$$
C_1 (1+r)e^{-(n-2)r}\le Q_n(r)\le C_2 (1+r)e^{-(n-2)r}
$$
so that  $Q_n$ is not constant in $r$. Actually, it is an analytic function of $r^2$ with a maximum at $r=0$. In three dimensions we have the simple expression
$$
Q_3(r)=\frac{C r}{\Sh r}.
$$
This function has a maximum at $r=0$ ($x=o$) and decreases rapidly as $r$ grows. This is in sharp contrast with the Euclidean case where we get local homogeneity at this level, i.e.,
$$
\lim_{t\to\infty}t^{n/2}E_t(x;0)= (4\pi)^{n/2}>0.
$$

\noindent (ii) In order to prove our negative result, we take as second solution the fundamental solution $u_1(x,t)=P_t(x,x_a)$  centered at $x_a=(a,0,\dots,0)$, and we compare it with $P_t(x,o)$. We  consider the points $x'=(r,0,...,0)$ in the geodesic line joining $P_0=o$ and $P_a=x_a$, and we get for some $r>a$
\begin{equation}\label{notconv.dm}
\lim_{t\to\infty}t^{3/2}e^{\lambda_1 t} |u_1(x',t) - P_t(x';o)|= Q(r-a)-Q(r)>0.
\end{equation}
Therefore, the $L^\infty$ version of Theorem \ref{thm.conv.hypb.pos} is not true.
%
%


\subsection{Study of the asymptotic mass error}

We have seen that there is a lack of uniform convergence near the origin that forces the expected $L^\infty$ convergence to fail. We want to prove a similar lack of sharp convergence for the $L^1$ norm.  Since we already know that the mass concentrated on bounded balls goes quickly to zero as $t\to\infty$, it is natural to enquire where and how the mass error is kept. The answer is that there is a non-vanishing  error on and around the so-called mass line $r=(n-1)t$. This is what we are going to show here. We work out the complete detail only in dimension $n=3$ for simplicity.

\medskip

{\bf Study of the solution to a displaced mass.}  Again, we take as solution $u$  the fundamental solution $u_1$ centered at $x_a=(a,0,\dots,0)$ and compare it with $P_t(x)=G_t(r)$, the fundamental solution centered at $x=o$. Our proofs will be based in finding bounds for the displaced solution compared with the fundamental solution along the axis of displacement. Note that the mass of the displaced solution will concentrate of a thin annulus at distance $(n-1)t$ from $x_a$ that will look like a distorted annulus as seen from the original center $x=o$.

\medskip

\noindent (i) {\bf  Control along the axis.} We consider the displaced solution with initial mass at $x_a$.  It value along the $x_1$ axis is given by
$$
u(x,t)= G_t(r-a)= Ct^{3/2}e^{-t}e^{-(r-a)^2/4t}\frac{r-a}{\Sh(r-a)},
$$
where $x=(r,0,\dots,0)$. Hence for $r>a$
$$
\frac{G_t(r-a)}{G_t(a)}=e^{\frac{ar}{2t}}e^{-\frac{a^2}{4t}}\,\frac{(r-a)\,\Sh r}{r\,\Sh(r-a)}.
$$
There is a bound for the last quotient in the interval $a<r<kt$
$$
\frac{(r-a)\,\Sh r}{r\,\Sh(r-a)}\le ce^a
$$
hence for all $a<r<kt$ we have
$$
\frac{G_t(r-a)}{G_t(a)}\le e^{\frac{ar}{2t}} ce^a=\le ce^{(1+k/2)a}.
$$
This includes the mass region $x-2t\sim \xi t^{1/2}$ where we have for large $t$
$$
\frac{G_t(r-a)}{G_t(a)}\sim e^{2a}>1.
$$
Note that there is a further correction factor $e^{-a^2/4t^2}$ that does not count for large $t\gg a$. Hence, in the interesting interval in the $x_1$-axis $G_t(r-a)$ is larger than and proportional to $G_t(r)$.

\medskip

\noindent (ii) \bf Control in a region. \rm We need to extend this result to a region around the interval in the $x_1$ axis in order to control the difference of mass of the two solutions. We will take an annulus with center $O$ and restrict it to the cone along the $x_1$ axis with  a small angle $\phi_0$ (amplitude). We call this region $R$.

In order to calculate the difference of solutions we have to compare for fixed radius $r\sim 2t$ the values of $u(x,t)=G_t(d(x_a,x))$ with $P_t(x)=G_t(r)$ for any point $x$ at distance $r$ from the origin $O$ and angle $0\le \phi\le \phi_0$. Let us call L the distance $d(x_a,x)$
We note that $L=L(a,\phi,r)$ is always located in the interval $r-a<L<r+a$. Indeed, for the calculation of quotients $u/G$ we need to estimate  $a(\phi)=r-L(a,\phi)$. For angle $\phi=0$ we have $L=r-a$ and the calculation of quotients $u/G$ has been just done and amounts to $ce^{2a}$ in the region $R$. For angle $\phi>0$ we need to estimate $a(\phi)$ from below.

In $n=2$ we may use the formula of hyperbolic geometry for triangles that reads
$$
\Ch L= \Ch r \,\Ch a -\Sh r \, \Sh a \,\cos \phi,
$$
We are interested in $a$ of small size and $r$ and $L$ very large. Then
$$
\Ch r \,\Ch a -\Sh r \, \Sh a \,\cos \phi=\frac14 e^{r+a}(1-\cos\phi)+ \frac14 e^{r-a}(1+\cos\phi)+o(1).
$$
We want $a(\phi)=r-L\ge a/2$ hence approximately
$$
e^{r+a}(1-\cos\phi)+ e^{r-a}(1+\cos\phi)<2e^{-a/2}e^r
$$
We need
$$
2e^{-a/2}\ge  e^a(1-\cos\phi)+ e^{-a}(1+\cos\phi)=2e^{-a}+ (e^a-e^{-a})(1-\cos\phi)
$$
hence
$$
e^{-a/2}(1-e^{-a/2})\ge \Sh a \,(1-\cos\phi).
$$
This is true for $\phi=0$ and it will be true for some $0<\phi<\phi_0$
that depends on $a$ but not on $r$ as long as $r$ is very large.

\medskip

\noindent (iii) \bf Mass estimate in the mass line region. \rm
The mass contained in that region is then proportional and larger for $u$ than for $G_t(r)$
by a constant factor.

We conclude that there exists a function $F(a)>0$ such that for all large $t\gg 1$
$$
\|G_t(r-a)-G_t(r)\|_1\ge F(a)
$$
which is the worst case for a bound of the form
$$
\|u(r,t)-G_t(r)\|_1\ge f(a,t)
$$
in the class of initial data which are bounded, nonnegative with compact support.

\medskip

\noindent (iv) {\bf Far field estimates. } On the other hand, in the far field we have the bound with respect to a fundamental solution that is a bit delayed in time:
$$
\frac{G_t(r-a)}{G_{t+\ve}(r)}=\left(\frac{t+\ve}{t}\right)^{3/2}\,\frac{(r-a)\Sh r}{r\Sh(r-a)}
\,e^{\frac{r^2}{4(t+\ve)}-\frac{(r-a)^2}{4t}}.
$$
But
$$
\frac{r^2}{4(t+\ve)}-\frac{(r-a)^2}{4t}=\frac1{4t(t+\ve)}\{-\ve r^2+ 2ra(t+\ve)-a^2(t+\ve) \}
$$
If $t$ is large, $\ve$ small and $r>kt$ we then have
$$
\frac{G_t(r-a)}{G_{t+\ve}(r)}\le 2e^a e^{-r(\ve r-2at)/4t^2}
$$
which is bounded uniformly if $r>(2a/\ve)t$. This gives a uniform bound and proves at the same time that the mass in that region is negligible for both solutions if $k$ is large.

From the analysis we get the following conclusion

\begin{theorem} \label{thm.nonconv.7} There is a positive constant $c=c(a)\in (0,1)$ such that for all large $t$
\begin{equation}
\displaystyle \|(u(\cdot,t)-G_t(\cdot))_+\|_{L^1(\mathbb{H}^3)}\ge c(a)>0.
\end{equation}
\end{theorem}

\subsection{Nonradial Counterexample. Second derivative in $a$.} We want to prove that there is no convergence to the fundamental solution centered at a certain center of mass, which was a possible remedy to the negative result. We consider the solution $u_1$  with initial data $1/2$ Dirac mass at $P_1=(1,0,\dots)$ and $1/2$ Dirac mass at $P_2=(-1,0,\dots)$ and prove that the negative result of Theorem \ref{thm.nonconv.7} still holds. Note that for large $t$ the solution mass will be concentrated on the union of two closely located distorted annuli.

To begin with,  at a point $r>1$ along the $x_1$ axis we get
$$
u_1(r,t)=\frac1{2}G_t(r-a)+\frac1{2}G_t(r+a)\,,
$$
so that for the region of mass interest $r\sim 2t$ we get
$$
u_1(r,t)\sim \frac1{2}e^{2a} G_t(r)+\frac1{2}e^{2a}G_t(r)=\Ch(2a)G_t(r)\,,
$$
which is larger than $G_t(r)$ by a factor  $K=\Ch(2a)-1>0$. Note that there is a further correction factor $e^{-a^2/4t^2}$ that does not count for large $t\gg a$.

We also have
$$
\log P(r,a)= \log C- (3/2) \log t+  \log (r/\Sh r) -t -\frac{r^2}{4t}.
$$
So that
$$
\frac{\partial_a P(r+a,t)}{P(r+a)}= \frac1{r+a}-\frac{\Ch (r+a)}{\Sh (r+a)}-\frac{r+a}{2t}.
$$
and also
$$
\frac{\partial^2_a P}{P}-\left(\frac{\partial_a P}{P}\right)^2=
-\frac1{(r+a)^2}-\frac{1}{\Sh^2 (r+a)}-\frac1{2t}.
$$
For the interesting region $r\sim 2t$ we have $\partial_a P/P\sim -2$ and $\partial^2_a P/P\sim
(\partial_a P/P)^2\sim 4$. This allows to produce the  counterexample.

\section{The complete equation}\label{sec.complete}

We now take a look at the complete equation
\begin{equation}
\partial_t u=\Delta_g u+ f\,,
\end{equation}
where the novelty is the forcing term $f=f(x,t)$. It is well known that, similarly to the Euclidean case, the complete equation in $\Hn$ generates an evolution process for all $f\in L^1(0,\infty; L^p(\Hn)$ with $1\le p<\infty$. Moreover, the following basic estimate is true for any two solutions $u_1,u_2$, with data $f_1,f_2$, and initial conditions $u_{01}, u_{02}.$
$$
\|u_1(t)-u_2(t)\|_{1}\le \|u_1(0)-u_2(0)\|_{1}+
\int_0^t \|f_1(s)-f_2(s)\|_{1}\,ds\,.
$$
with norms in ${L^1(\Hn)}$. In a rather standard way we can combine this estimate with Theorem \ref{thm.conv.hypb.pos} and get

\begin{theorem} \label{thm.conv.hypb.compl}  Let $u(x,t)$ be a solution of the complete heat equation in the hyperbolic space whose initial function $u_0\in L^1(\mathbb{H}^n)$ is radially symmetric in geodesic coordinates around a pole $o\in \Hn$ and has mass $M_0$, and let $f=f(r,t)$ belong to $L^1_{t,x}(\Hn\times (0,\infty))$. Then we have
\begin{equation}
\displaystyle \|u(x,t)-MP_t(x)\|_{L^1(\mathbb{H}^n)}\to 0 \quad \mbox{as \ $t\to\infty$}\,,
\end{equation}
where $M$ is the accumulated mass,
\begin{equation}
M=M_0+\int_0^\infty\int_{\Hn} f(r,t)\,d\mu\,dt.
\end{equation}
\end{theorem}

Results about convergence in sup norm need extra assumptions on the decay of $f$ for large time. We skip further details.

\section{Comments and problems}\label{sec.comm}

\noindent $\bullet$ If we compare the mode of asymptotic convergence in the 3 main geometries, we find quite different types of convergence: (I) convergence with a polynomial rate to the Gaussian profile in the Euclidean case (curvature $K=0$), (II) convergence to the constant state with exponential rate in the spherical case (curvature $K>0$), and (III) the kind of convergence we have described here in the hyperbolic geometry (curvature $K<0$). To note also that the heat flow in bounded domains of $\ren$ looks like the $K<0$ case, in the sense for both Dirichlet and Neumann problems convergence to the asymptotic profile happens with an exponential rate.

\noindent $\bullet$  We pose as an open problem finding a clear proof of Theorem \ref{thm.conv.hypb} in dimensions $n\ge 2$ other than 3.

\noindent $\bullet$ The investigation of the precise large time behaviour of nonlinear heat flows on hyperbolic space has been done by the author in the case of the Porous Medium Equation in \cite{Vazhyp15}, with quite different results. Thus, no equivalent to mass escape is found, either in the radial or in the horospheric class of solutions. Also, the $p$-Laplacian equation was considered in \cite{Vazhyp15}. Radial fast diffusion on the hyperbolic space was considered in \cite{GM14}. To note that for flows in the Euclidean space the nonlinear flows corresponding to the PME  and Fast Diffusion Equation are not so much at variance with the linear flow (the heat equation in $\ren$), see \cite{VazBook, VazCIME}.

On the same vein, there is a stark contrast when comparing the horospheric solutions in the linear heat equation case (see Section \ref{sec.horo}) and the nonlinear PME analysis done in \cite{Vazhyp15}. We recall that in this setting there are many explicit formulas that help understand the evolution.

\noindent $\bullet$ In the Porous Medium Equation case the analysis has been extended to flows
in other manifolds with negative curvature of the type called  model manifolds, under conditions on the behaviour of the negative curvature ``at infinity'', see \cite{GMV17, GMV18}. Such study is to be done in the heat equation case.

\noindent $\bullet$ It is interesting to insert the quantitative effect of a negative curvature $K<0$ and examine the limit $K\to 0$ to see how the hyperbolic effects disappear. The heat equation becomes
$$
\partial_t u= \Delta_g u \equiv (\sinh (r/R))^{1-n}\frac{\partial}{\partial r}
\left((\sinh (r/R))^{n-1}\frac{\partial u}{\partial r}\right)
+\frac{1}{R^2\,(\sinh (r/R))^2} \Delta_{\mathbb{S}^{n-1}}u\,,
$$
where $R^2=-1/K$.

\noindent $\bullet$ It will be interesting to better understand the local space inhomogeneity of the fundamental solutions, that for $n=3$ reads $Q(r)=r/\Sh r$. As we saw, it plays a role in the negative convergence result.

\noindent $\bullet$
A  natural generalization would be to consider variable curvature manifolds such Cartan-Hadamard manifolds. It has been considered by the author and collaborators in previous papers. We mention in this section a natural framework for similar results: the case of {\sl any} Riemannian symmetric space of non-compact type. When dealing with real hyperbolic spaces, it is natural to consider actually the larger class of all rank 1 symmetric spaces of non-compact type. The table below provides all of them.
\vspace{2mm}

\centerline{\begin{tabular}{|c|c|c|c|c|}
\hline
\hspace{10mm}
&\,$\mathbb{H}^{\ssf d}\!=\msb\text{\rm H}^{\ssf d}(\re)$
&\hspace{5mm}$\text{\rm H}^{\ssf d}(\C)$\hspace{5mm}
&\hspace{5mm}$\text{\rm H}^{\ssf d}(\H)$\hspace{5mm}
&\hspace{5mm}$\text{\rm H}^{\ssf 2}(\O)$\hspace{5mm}$\vphantom{\frac||}$\\
\hline
$n$&$d$&$2\ssf d$&$4\ssf d$&$16\vphantom{\frac||}$
\\\hline
$\rho$&$\frac{d\vsf-1}2$&$d$&$2\ssf d\msb+\!1$&$11\vphantom{\frac||}$
\\\hline
\end{tabular}}
\vspace{2mm}
In the table, we denote by $n$ the Riemannian dimension and by $\rho$, the bottom of the spectrum of the Laplace operator. The idea behind is that (after Helgason, see for instance \cite{Helgason1984}) the spherical analysis on all these spaces is completely similar and thus provides similar bounds on the heat kernel. More precisely, we refer to the work of  Anker and Ostellari \cite{AnkOst} where precise estimates for the heat kernel on symmetric spaces of non-compact of any rank (see their Main Theorem). Furthermore, the result on heat propagation concentrated on a precise annulus in Section 4 of the same work is consistent with ours. It has to be noticed that such results are also true on Damek-Ricci spaces (as discovered in \cite{ankerdamek}). Extending our asymptotic analysis to those settings is therefore a natural task.

\noindent $\bullet$ Different variants of the basic  heat equation have been studied. Let us give an example of a heat equation with a reaction term: the equation $u_t=\Delta_g(u)+ f(u)$ has been studied in \cite{MatPunTes} in the KPP case $f(u)=u(1-u)$, and the presence of travelling waves as asymptotic profiles is established. This is quite different from the behaviour described here in the free case $f=0$.

\vskip 1cm

\noindent {\textbf{\large \sc Acknowledgments.}} Partially funded by Project  MTM2014-52240-P (Spain). Partially performed as a Visiting Professor at Univ. Complutense de Madrid. The observations about non-compact rank 1 symmetric spaces are largely due to Y. Sire whom I thank so much for his attention to this work. I also thank G. Grillo for interesting conversations.

\

\bibliographystyle{amsplain} 

\

\noindent {\sc Address:}

\noindent Juan Luis V\'azquez. Departamento de Matem\'{a}ticas, Universidad
Aut\'{o}noma de Madrid,\\ Campus de Cantoblanco, 28049 Madrid, Spain.  e-mail address:~\texttt{juanluis.vazquez@uam.es}

\


\end{document}